\documentclass[12pt,a4paper]{article}

\usepackage[british]{babel}
\usepackage{graphicx}
\usepackage{amsmath}
\usepackage[]{mathtools}
\usepackage{geometry}
\usepackage{amssymb}
\usepackage{epstopdf}
\usepackage{enumitem}
\usepackage{xcolor}

\begin{document}

\title{Using a library of chemical reactions to fit systems of ordinary differential equations to agent-based models: a machine learning approach}
\author{Pamela M. Burrage\footnotemark[1], \and Hasitha N. Weerasinghe{\footnotemark[1]}, \and Kevin Burrage{\footnotemark[1] \footnotemark[2]}}

\maketitle
\renewcommand{\thefootnote}{\fnsymbol{footnote}}

   \footnotetext[1]{School of Mathematical Sciences,
      Queensland University of Technology (QUT), Australia. kevin.burrage@qut.edu.au}
      \footnotetext[2]{Visiting Professor, Department of Computer Science, University of Oxford, United Kingdom}

   \renewcommand{\thefootnote}{\arabic{footnote}}

\begin{abstract}

In this paper we introduce a new method based on a library of chemical reactions for constructing a system of ordinary differential equations from stochastic simulations arising from an agent-based model. The advantage of this approach is that this library respects any coupling between systems components, whereas the SINDy algorithm (introduced by Brunton, Proctor and Kutz) treats the individual components as decoupled from one another. Another advantage of our approach is that we can use a non-negative least squares algorithm to find the non-negative rate constants in a very robust, stable and simple manner. We illustrate our ideas on an agent-based model of tumour growth on a 2D lattice.

\end{abstract}




\section{Introduction}

As Brunton and Kutz \cite{BruKu} state, ``Modern dynamical systems are currently undergoing a renaissance with analytical derivatives and first principles models giving way to data-driven approaches. The confluence of big data and machine learning is driving a paradigm shift in the analysis and understanding of dynamical systems.'' In this paper, we will use the Sparse Identification of Nonlinear Dynamics (SINDy) \cite{BruPrKu} to analyse an agent-based spatial model of tumour growth, The general aim is to discover dynamical systems models from data via automated discovery of the underlying governing equations. The fundamental idea underpinning SINDy is that many dynamical systems have an underlying Ordinary Differential Equation (ODE) form 
\begin{equation}
y'(t) = f(y(t)), \quad f,\> y \in  \mathbb{R}^d ,
\label{eq:eq1}
\end{equation}
where the right hand side  has just a few active terms (often linear and quadratic) and we want to represent $f$ as a linear combination of these terms (library of candidate functions).

The formulation in (\ref{eq:eq1}) is through a column representation so $y = (y_1,\cdots,y_d)^\top$. However, it is sometimes more convenient in the SINDy representation to use a row-wise formulation so, for example,
$$ y' = (y_1',\cdots,y_d'), \quad f(y) = (f_1(y),\cdots,f_d(y)) \in \mathbb{R}^d. $$
We now write the library of $D$ scalar candidate functions ($\theta_j(y), \> j=1,\cdots,D$) used to approximate the right hand side of (\ref{eq:eq1}) as
$$ \theta(y) = (\theta_1(y),\cdots,\theta_D(y)) \in \mathbb{R}^D.$$
Let $\xi_j \neq 0, \> j=1,\cdots,d $ represent the coefficients of the $D$ candidate functions. We  let $\Sigma$ be the matrix
$$ \Sigma = (\xi_1,\cdots,\xi_d) \in \mathbb{R}^{D \times d}.$$
Note that individual components of the $\xi_j$ can be different from 1 and represent a general coefficient. In the row-wise formulation this leads to
\begin{equation}
y'  \approx \theta(y) \Sigma.
\label{eq:eq2}
\end{equation}

The aim here is to have $D$ as small as possible consistent with the data and for the $\xi_j$ to have just a few non-zero active terms. In order to find these terms we can use time series data collected from the numerical solution of (\ref{eq:eq1}) at a set of $N$ time points $t_1,\cdots,t_N$, or indeed collected from say an agent-based model (ABM). Thus an $N \times d$ data matrix can be represented by  $Y$,
$$ Y = (y(t_1), \cdots,y(t_N))^\top \> \in \mathbb{R}^{N \times d},$$
along with an $N \times d$ matrix of numerical derivative approximations
$$ Y' = (y'(t_1),\cdots,y'(t_N))^\top \> \in \mathbb{R}^{N \times d}.$$
We now construct $\theta(Y) \in \mathbb{R}^{N \times D}$, a function of the data matrix $Y$, to represent
\begin{equation}
Y' = \theta(Y) \Sigma.
\label{eq:eq3}
\end{equation}

This approach has also been extended to the identification of partial differential equations \cite{RudBrPrKu}. In this case a partial differential equation, in say one spatial dimension, can be expressed in the library of $D$ candidate terms as
\begin{equation}
\hat{Y}_t = \theta(\hat{Y}) \Sigma
\label{eq:eq4}
\end{equation}
where  $\theta(\hat{Y}) \in \mathbb{R}^{mN \times D}$ can not only represent functions of $\hat{Y}$ but also various partial derivatives of $\hat{Y}$.

Here $\hat{Y}$ represents data collected over $N$ time points and $m$ spatial locations. As before,  the vectors $\xi_j$ should be sparse. Derivatives are again approximated using finite differences.  Both formulations (\ref{eq:eq3}) and (\ref{eq:eq4}) are regression problems and may be ill-conditioned if some of the columns of $\theta(Y)$ are close to being parallel or linearly dependent upon one another. Of course we can choose our library elements to avoid this linear dependence if necessary. This ill-conditioning can be accentuated by, for example, approximation errors in the derivatives. 

At its simplest level, the regression problem can be solved by linear least squares leading to the solution of the normal equations.  Thus, given the rectangular system $Ax = b$, the linear least square solution minimises $||b-Ax||^2_2$ through the normal equations 
$$ (A^\top A) x = A^\top b$$
i.e.
\begin{equation}
x = (A^\top A)^{-1} A^\top b.
\label{eq:eq6}
\end{equation}
The \texttt{lsqr} algorithm in MATLAB is an adaptation of the Conjugate Gradient method for rectangular systems that is generally more reliable than computing (\ref{eq:eq6}) \cite{BarBeCh}. However, the standard least squares approach can still be very sensitive to random errors in approximating derivatives (encoded here in the $b$ vector). Furthermore, if there are correlations between some of the columns of the $A$ matrix, then $A^\top A$ can be close to being singular. Ridge regression \cite{HoeKe} attempts to overcome this issue of multi-collinearity by introducing a positive ridge parameter $k$ so that an estimate of the regression coefficients is given by 
\begin{equation}
x = (A^\top A + kI)^{-1} A^\top b.
\label{eq:eq7}
\end{equation}
This can result in smaller mean squared errors compared to the solution in (\ref{eq:eq6}). In addition, it is possible to use a ridge regularisation and a ridge regression parameter $\lambda$ so that the minimisation problem becomes $\arg \min(||Ax-b||^2_2 + \lambda ||x||^2_2)$. This can then be solved by the Lasso algorithm \cite{Tibshirani}.
Brunton and Kutz \cite{BruKu} discuss further refinements of the SINDy problem by asking how to find sparse vectors   with a small solution residual. They and Rudy et al. \cite{RudBrPrKu} propose ridge regression with hard thresholding. Thus, given a tolerance and a threshold  they iteratively refine ridge-regression by removing values in the active components that are less than the tolerance in absolute value. 


SINDy has been applied in a wide variety of applications. These include the identification of high dimension systems in fluid flows \cite{LoiBr}, the construction of reduced order methods in plasma dynamics \cite{KapMoHaBr},  nonlinear optics \cite{SorSyTu}, nonlinear stochastic modelling \cite{LibWoMa} and thermal convection \cite{Loiseau}. In the fields of modelling biological and chemical systems, Brunton and Kutz \cite{BruKu} note that many dynamical systems are described by rational nonlinear functions - via Hill functions. In order to apply SINDy in this setting, it is necessary to use implicit differential equations and the library of candidate functions is generalised to include state and derivative dependencies. Finally, other approaches for discovering dynamical systems have been pioneered including equation-free modelling \cite{KevGeHy}. 

Nardini et al. \cite{NarBaSiFl} give a nice tutorial on learning differential equation models from agent-based model simulations. They focus on a one-dimensional birth-death model and a two-dimensional SIR model. In the first case most regression approaches will work robustly and they infer a quadratic or order 4 right hand side. However, in the latter case, inferring two active components with a library $\theta(S,I) = (S, S^2, I, I^2,SI)$ is not trivial nor robust. An iterative pruning method based on FISTA \cite{FISTA} is used in conjunction with Lasso and there is considerable delicacy in selecting the regularisation parameter $\lambda$ to be used within Lasso. We will discuss these issues in more detail in section 2. 

We will apply the SINDy approach to an agent-based simulation model of cancerous tumour growth in which we model the interaction between tumour cells and normal healthy cells with the aim of deriving a coupled system of ODEs from the time series data arising from the agent-based model. We simulate the agent-based model 500 times and average the output at each time point in order to obtain the data in the appropriate form. We explore the behaviour in terms of the number of time points.

In section 2, we describe a modification of the SINDy approach whereby we construct a library of chemical reactions that naturally couples the resulting system components of the ensuing ODE.  The standard SINDy approach does not do this. In section 3 we describe the agent-based model in greater detail \cite{Hthesis}. Section 4 gives numerical results and discussion while  section 5 offers some  conclusions.

\section{A library of chemical reactions}

One issue with the standard SINDy approach outlined in Brunton and Kutz \cite{BruKu} and Nardini et al. \cite{NarBaSiFl} is that the differential system constructed is done equation by equation and so has its system components decoupled from one another. If an agent-based model is constructed based on a set of chemical reactions and the same library of functions is used for each component, for example, then, in general, the stoichiometric vectors couple the components.  This is not reflected in the standard SINDy approach, nor is it possible to infer the underlying chemical reactions that underpin the agent-based simulations. 

In order to address this issue, we move away from a library of functions that generally are not coupled between the system components and work with a library of chemical reactions. For example, in the case of two components $x$ and $z$, the complete set of unimolecular and bimolecular reactions are as follows. Note we do not consider tri-molecular  and higher order reactions.

\vspace{1cm}
\noindent
\textbf{Unimolecular reactions}
$$ X \rightarrow 0, \> Z \rightarrow 0, \> X \rightarrow Z,  \>Z \rightarrow X, \> X \rightarrow X+X, $$
$$  X \rightarrow X+Z, \> Z \rightarrow Z+Z, \> Z \rightarrow Z+X $$
with stoichiometric vectors 
$$ [-1,0]^\top, [0,-1]^\top, [-1,1]^\top, [1,-1]^\top, [1,0]^\top, [0,1]^\top, [0,1]^\top, [1,0]^\top $$
and propensity functions (ignoring rate constants)
$$ x, \> z, \>  x, \> z, \>  x, \> x, \>  z, \> z. $$

\vspace{1cm}
\noindent
\textbf{Bimolecular reactions}
$$ X+X \rightarrow 0, \> Z+Z \rightarrow 0, \> X+X \rightarrow X, \> X+X \rightarrow Z, \> Z+Z \rightarrow Z, $$
$$ Z+Z \rightarrow X, \> X+Z \rightarrow 0, \> X+Z \rightarrow X, \> X+Z \rightarrow Z $$
with stoichiometric vectors
$$ [-2,0]^\top, [0,-2]^\top, [-1,0]^\top, [-2,1]^\top, [0,-1]^\top, [1,-2]^\top, [-1,-1]^\top,  [0,-1]^\top, [-1,0]^\top $$
and propensity functions
$$ \frac{1}{2} x^2, \> \frac{1}{2} z^2, \> \frac{1}{2} x^2, \> \frac{1}{2} x^2, \> \frac{1}{2} z^2, \> \frac{1}{2} z^2, \> xz, \> xz, \> xz. $$
Note that the propensity functions are equal to the products of the reactants on the left hand side of the reaction times the associated rate constant. For example, the propensity function associated with the reaction $X+Z \rightarrow Z$ is $xz$. If the associated rate constant is $k$, then the propensity function is $kxz$.
Thus the library of 17 possible reactions is given by
$$ \theta(x,z) = \left( \left(\begin{array}{r} -x \\ 0 \end{array}\right), \> \left(\begin{array}{r} 0 \\ -z \end{array}\right) \> \left(\begin{array}{r} -x \\ x \end{array}\right), \> \left(\begin{array}{r} z \\ -z \end{array}\right), \> \left(\begin{array}{r} x \\ 0 \end{array}\right), \> \left(\begin{array}{r} 0 \\ x \end{array}\right), \>  \left(\begin{array}{r} 0 \\ z \end{array}\right), \right. $$ 
$$  \quad  \quad \quad \quad \quad \quad  \quad  \left.  \left(\begin{array}{r} z \\ 0 \end{array}\right),  \> \left(\begin{array}{r} -x^2 \\ 0 \end{array}\right), \> \left(\begin{array}{r} 0 \\ -z^2 \end{array}\right), \> \left(\begin{array}{r} -\frac{1}{2}x^2 \\ 0 \end{array}\right), \> \left(\begin{array}{r} -x^2 \\ \frac{1}{2}x^2 \end{array}\right), \> \left(\begin{array}{r} 0 \\ -\frac{1}{2}z^2 \end{array}\right), \right. $$
$$ \quad   \quad  \left.  \left(\begin{array}{r} \frac{1}{2}z^2 \\ -z^2 \end{array}\right),  \>  \left(\begin{array}{r} -xz \\ -xz \end{array}\right), \>  \left(\begin{array}{r} 0 \\ -xz \end{array}\right), \>  \left(\begin{array}{r} -xz \\ 0 \end{array}\right) \right). $$
Under this library, the two-dimensional system of ODEs is 
\begin{eqnarray*}
x' &=& (k_5 -k_1 -k_3)x + (k_4+k_8)z - (k_9 + \frac{1}{2}k_{11} + k_{12})x^2 + \frac{1}{2} k_{14} z^2 \\
&& - (k_{15}+k_{17}) xz \\
z' &=& (k_7-k_2-k_4)z + (k_3+k_6)x - (k_{10}+\frac{1}{2}k_{13}+k_{14})z^2 + \frac{1}{2}k_{12}x^2 \\
&&  - (k_{15}+k_{16})xz.
\end{eqnarray*}
In general the system of ODEs that describes the dynamics of a set of $m$ reactions with stoichiometric vectors $\nu_j$ and propensity functions $a_j(X(t))$, where $X(t)$ is the state vector, is given by 
$$ X' = \sum_{j=1}^m \nu_j a_j(X).$$ 
Note that we have used the notation $k_j$ for the rate constants rather than the $\xi_j$. We can also see that there is a natural coupling between the two components via reactions 3, 4, 12, 14, 15. 
The least squares problem becomes 
$$ \min \left( \left| \left| \left( \begin{array}{c} X' \\ Z' \end{array} \right)^\top - \theta (X,Z) K \right| \right|^2_2 \right) $$
where $X$ and $Z$ are the data vectors and the vector of rate constants is  $K = (k_1,k_2,\cdots,k_{17})^\top.$ Here the vector $K \geq 0$ and so the regression problem has the constraint that all the $k_i$ should be nonnegative. This can be solved by a non-negative least squares algorithm, for example \texttt{lsqnonneg} in MATLAB, which is based on an iterative algorithm described in \cite{LawHa}.

Of course in many cases not all the reactions given in the library will appear in the model and so some of the $k_j$ will be 0. We discuss this along with numerical results in section 4.

\section{The tumour agent-based model}
Cancer is a disease that involves abnormal cell growth. The group of abnormal cells is called a tumour. If tumour cells spread to surrounding tissues or other organs then the disease can be life-threatening. The environment surrounded by a tumour is called the Tumour Micro-Environment (TME). The Extracellular Matrix (ECM) consists of proteins and provides cellular structural support by making the basement membrane and interstitial matrices. As can be seen in the cartoon in Figure 1,
there can be different types of cancer-associated cells and there may be T-cells present (innate or medication-induced) that are programmed to attack and kill cancerous cells.

\begin{figure}
\centering
\includegraphics[width=0.55\textwidth]{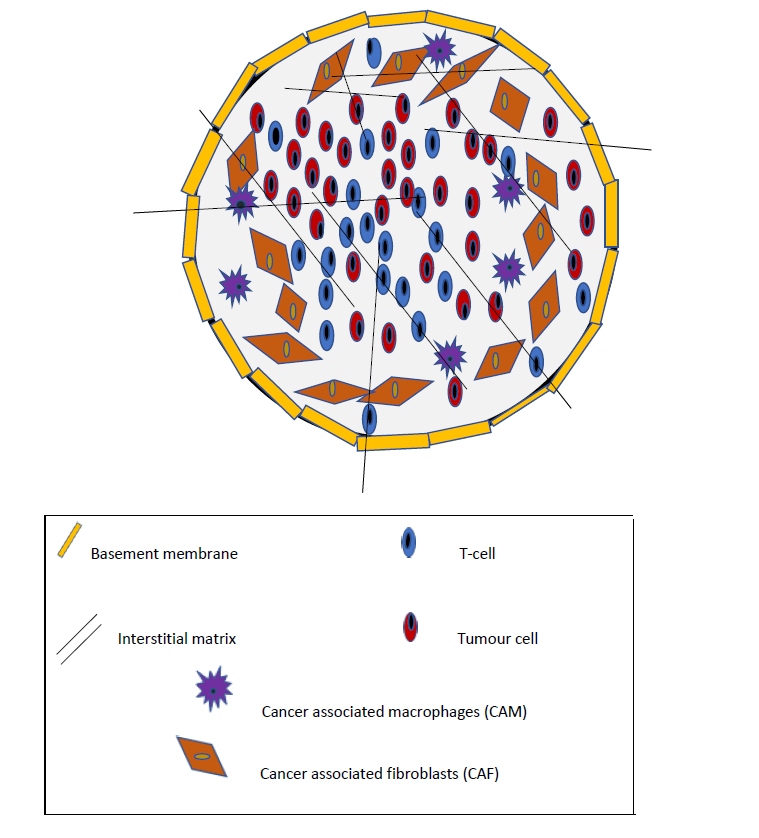}
\label{fig:figTME}
\caption{Cartoon showing the tumour microenvironment}
\end{figure}

In modelling cancer growth, even in very simple settings, there are key features that a model should represent \cite{HanWe2000}, \cite{HanWe2011}. In no particular order these include:
\begin{itemize}
\item tumour cells should have high proliferation rates in order to evade growth suppressors;
\item tumour cells should be able to proliferate uncontrollably in order to sustain proliferative signalling and allowing for replicative immortality;
\item tumour cells should have a low death rate;
\item tumour cells should have a high metabolism;
\item tumour cells should be able to diffuse to the boundary of the TME -- hence inducing angiogenesis, that is the formation of new blood vessels that can nourish the tumour;
\item and finally, tumour cells should be able to escape the TME thus giving rise to metastasis.
\end{itemize}

Our agent-based model consists of a square lattice with a $100 \times 100$ equidistant grid. It has two regions, an outer region and an inner region. At time zero, the boundary is fully fenced and the domain consists of healthy cells and immobile ECM proteins that represent obstacles. In Figure 2
the obstacles are marked by crosses and the healthy cells by blue circles. A single random cell mutates and the daughter cell becomes a tumour cell, marked as red.

\begin{figure}
\centering
\includegraphics[width=0.55\textwidth]{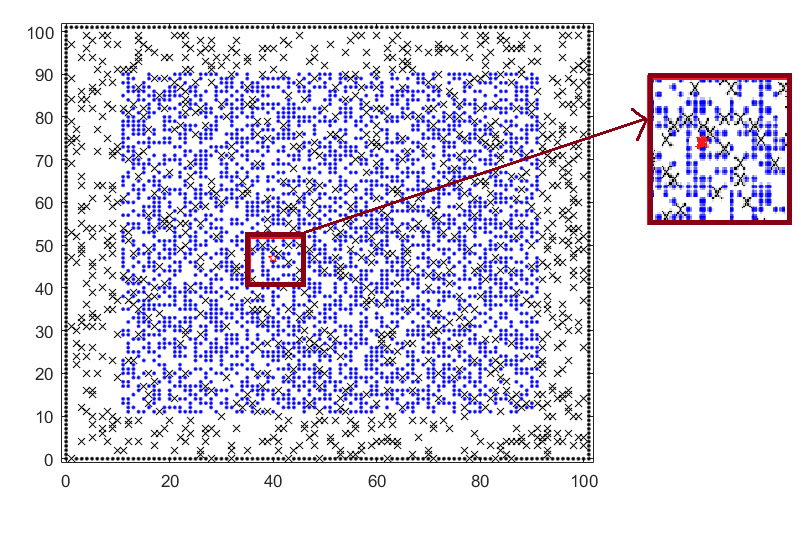}
\label{fig:figABM}
\caption{Simulation set up for the agent-based model}
\end{figure}

Every cell has a stickiness, jump radius, division age, competition rate, and probabilities of death, movement and ECM breakdown. When tumour cells reach the boundary of the lattice, they attach to the actin fence surrounding the lattice and secrete proteases to breakdown these components. The maximum number of cell divisions of a healthy cell is a parameter $K$. Tumour cells can divide more than $K$ times. A daughter cell can move to an empty grid position in the neighbourhood of a mother cell, it can also move to the grid position that a weaker cell may have occupied. Cells can also undergo jumps of different lengths. In summary, our model investigates tumour cell dynamics for various initial healthy cell densities, ECM protein densities and ECM breakdown probabilities. 

In Figure 3,
we show typical outputs from our agent-based model. The top row  shows the effects of different initial healthy cell densities (uniformly generated) and different obstacle densities  on the cell stickiness as a function of time, $t$, for $t \in [0,20]$. The second row shows the effects of these densities (healthy, ECM) on the jump radius as a function of time.  We give the tumour and healthy cell densities as a function of time for the agent-based model in section 4 where we compare them with the generated ODE.

\begin{figure}
\centering
$$\begin{array}{cc}
\includegraphics[width=0.45\textwidth]{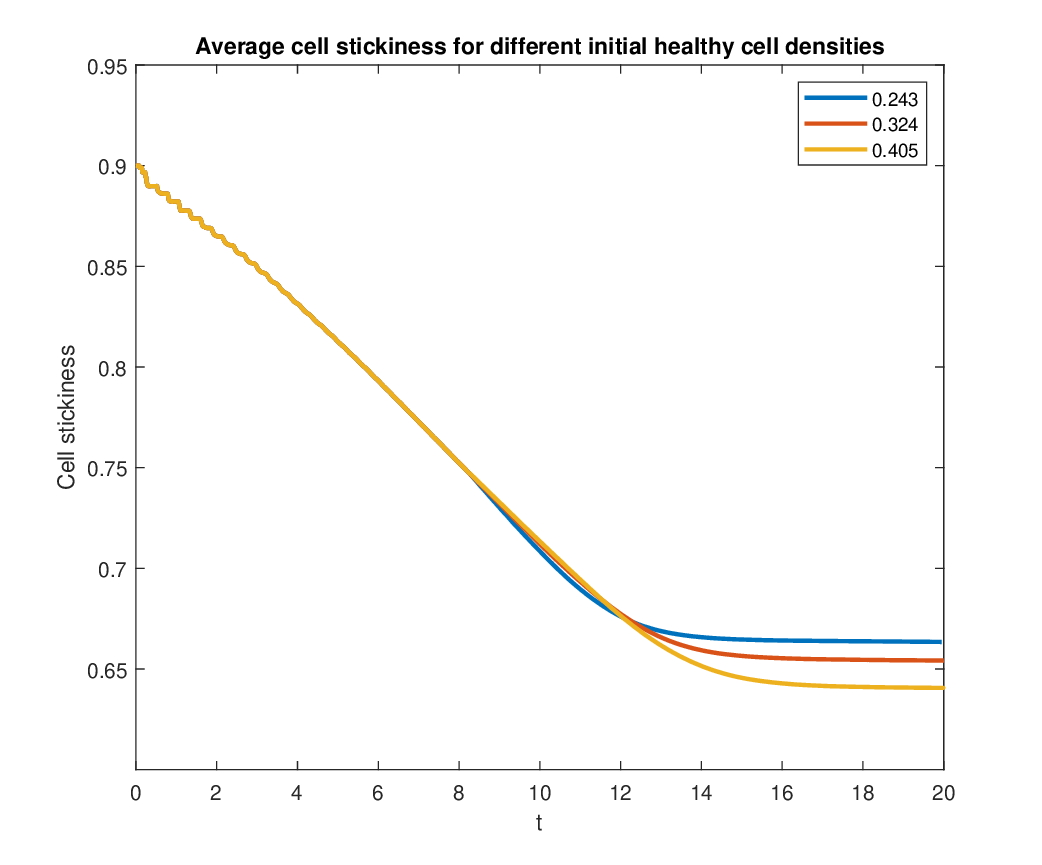} & \includegraphics[width=0.45\textwidth]{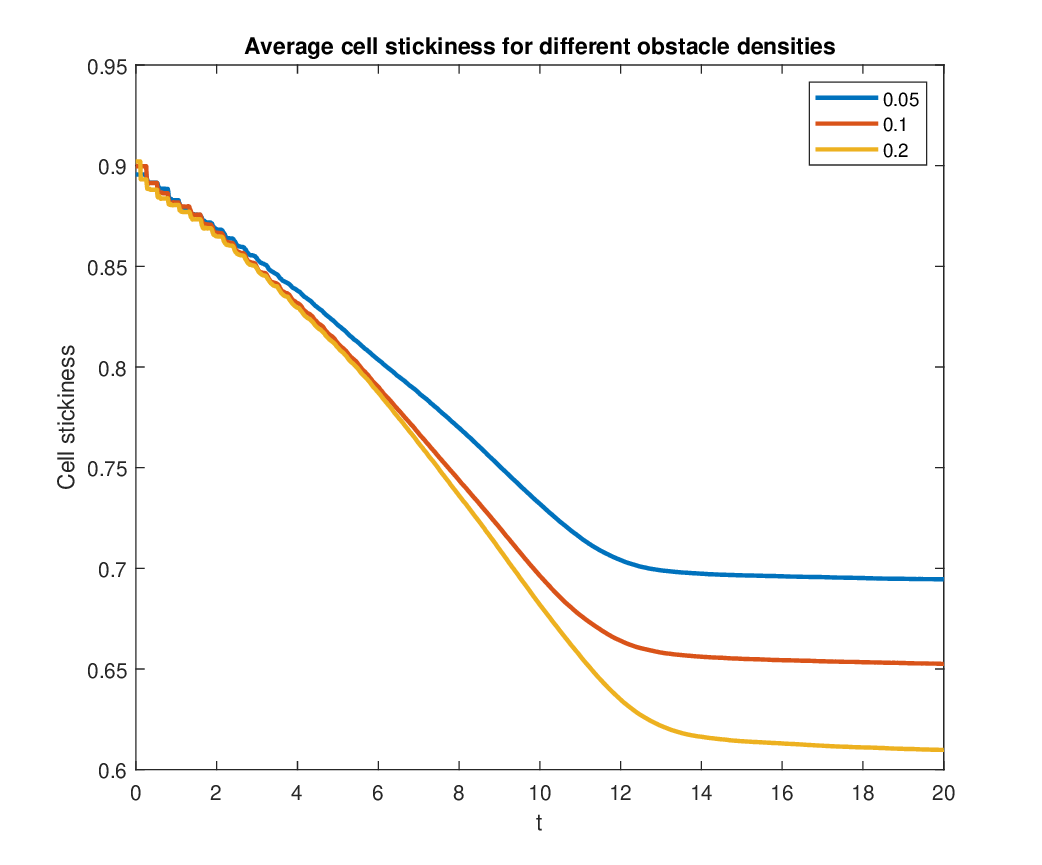}  \\
\includegraphics[width=0.45\textwidth]{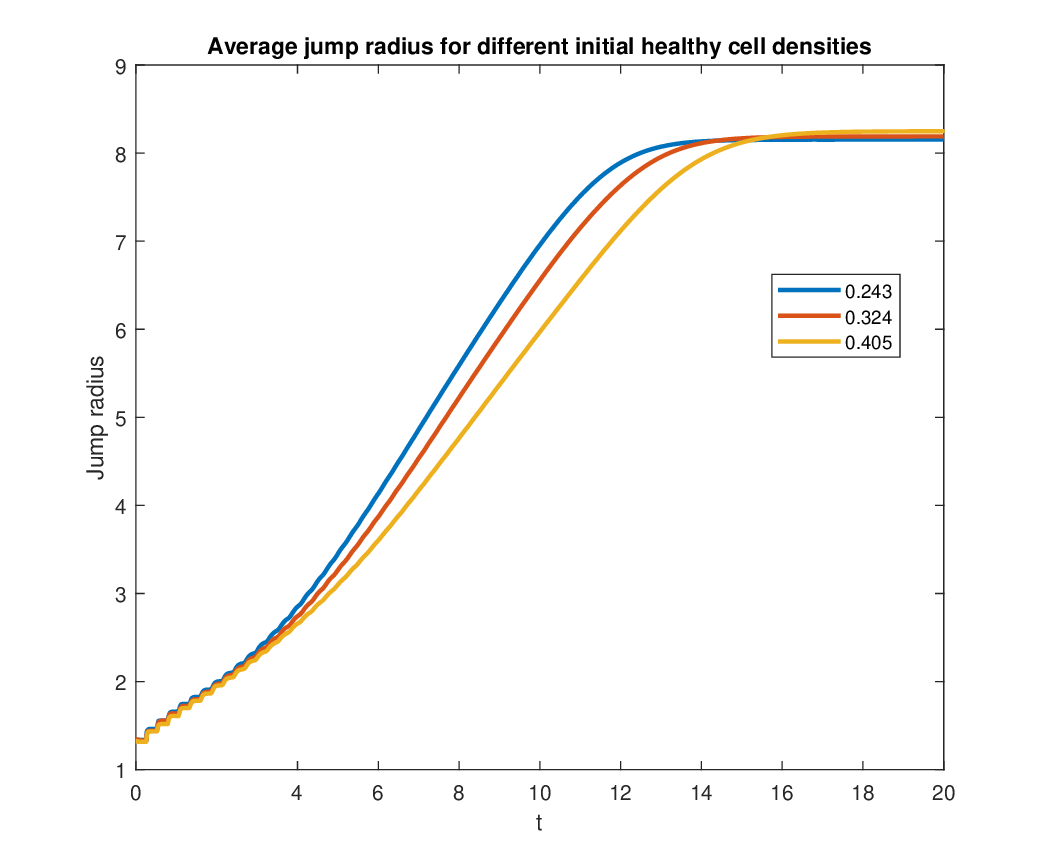} & \includegraphics[width=0.45\textwidth]{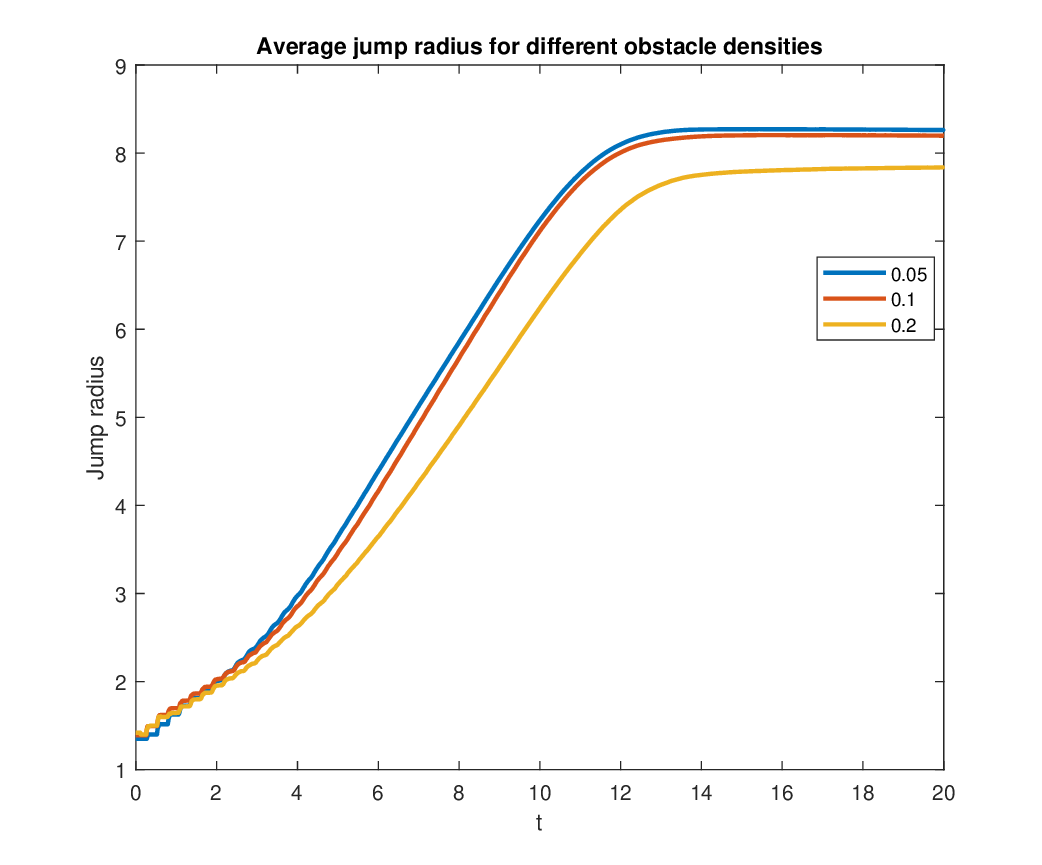}  
\end{array}
$$
\label{fig:figABMoutputs}
\caption{Sample outputs from the agent-based model in terms of effects on cell stickiness (top) and jump radius (bottom) in terms of initial healthy cell and ECM protein densities}
\end{figure}

\section{Implementation and numerical results}

We simulate the agent-based model and average the output at each time point. We take a time step size of $\frac{2}{135}$ and a time interval of $ [0, 26\frac{2}{3}]$, so that there are 1800 steps. We denote the number of tumour cells by $x$ and the number of healthy cells by $z$ and let $y = (x,z)^\top$. The initial healthy cell density was $z(0) = 0.324$, and the initial cancer cell density was $x(0) = 0.001$, corresponding to one cancer cell in the domain; the ECM protein density was 0.1 and the ECM breakdown probability was 0.5. The initial condition for healthy cells was chosen to give a nice balance between empty space and an initial crowding of healthy cells. A similar strategy to the initial healthy cell density was used for selecting the ECM protein density, while an ECM breakdown probability of 0.5 removes bias.

In the standard (non-coupling) approach, we choose our library of possible library functions $\theta(y)$ as
$$ \theta(y) = (x, x^2, z, z^2, xz).$$
In other words, we want to derive a system of ODEs in the form
\begin{eqnarray*}
x'(t) &=& a_1 x + a_2 x^2 + a_3 z + a_4 z^2 + a_5 xz \\
z'(t) &=& b_1 x + b_2 x^2 + b_3 z + b_4 z^2 + b_5 xz .
\end{eqnarray*}

We first form the data vectors by averaging the output at $t_1,\cdots,t_N$ over $k = 500$ simulations
\begin{eqnarray*}
\bar{X}(t) &=& \frac{1}{k} \left( \sum_{j=1}^k X_j(t_1),\cdots,\sum_{j=1}^k X_j(t_N) \right)^\top \\
\bar{Z}(t) &=& \frac{1}{k} \left( \sum_{j=1}^k Z_j(t_1),\cdots,\sum_{j=1}^k Z_j(t_N) \right)^\top .
\end{eqnarray*}
Then we form the numerical derivatives using finite differencing to obtain vectors $X', Z'$ of dimension as described in section 1. Our data model will have the form
\begin{eqnarray*}
X' &=& \theta(X,Z) \xi_x \\
Z' &=& \theta(X,Z) \xi_z.
\end{eqnarray*}
Since $\theta(X,Z)$ was evaluated at the $N$ time points $t_1,\cdots,t_N$ then 
$$ \theta(X,Z) = \left[ \begin{array}{ccccc}
\bar{X}(t_1) & \bar{X}^2(t_1) & \bar{Z}(t_1) & \bar{Z}^2(t_1) & \bar{X}(t_1)\bar{Z}(t_1) \\
\vdots & \vdots & \vdots & \vdots & \vdots \\
\bar{X}(t_N) & \bar{X}^2(t_N) & \bar{Z}(t_N) & \bar{Z}^2(t_N) & \bar{X}(t_N)\bar{Z}(t_N) 
\end{array} \right] . $$

In the second approach, we use the library of chemical reactions  that couples the $X$ and $Z$ components, and so we do not form two independent equations but one coupled equation.  The data model is 
$$ \left( \begin{array}{c} X \\ Z \end{array} \right)' = \theta_C (X,Z) \xi. $$
Here $ \left( \begin{array}{c} X \\ Z \end{array} \right)'$ is a $2N \times 1$ vector of derivative approximations, $\theta_C(X,Z)$ is a $2N \times 17$ matrix and $\xi$ is the $17 \times 1$ vector corresponding to the 17 non-negative rate constants associated with the 17 reactions.

We used a variety of methods, including Least Squares Regression, Ridge-Regression and Lasso, to determine the rate constants in the uncoupled case. However, these methods were nearly always unable to determine a stable system of ODEs. The only partial success was with the MATLAB code \texttt{lsqr}. In this case an  error message was returned  that convergence was not attained within the specified number of maximum iterations. But nevertheless, a solution was returned:
\begin{eqnarray}
x' &=& 2.5877x - 2.6283x^2 - 0.7384z + 2.3198 z^2 - 5.9795 xz  \nonumber \\
z' &=& 0.1407z - 0.4710 z^2 - 0.6011x + 0.6204 x^2 + 1.2618 xz.
\label{eq:eq1star}
\end{eqnarray}
The mean-squared error at $t = 26\frac{2}{3}$ was
$$ (x,z) = (1.8768(-04), 1.6369(-05))$$
based on 1800 time points. We note that we did not try the iterative FISTA algorithm in conjunction with Lasso. Rather we were trying to ascertain whether a simple numerical approach would yield meaningful results. 
See Figure 4.

\begin{figure}
\centering
\includegraphics[trim=10mm 70mm 20mm 70mm,clip,width=0.5\textwidth]{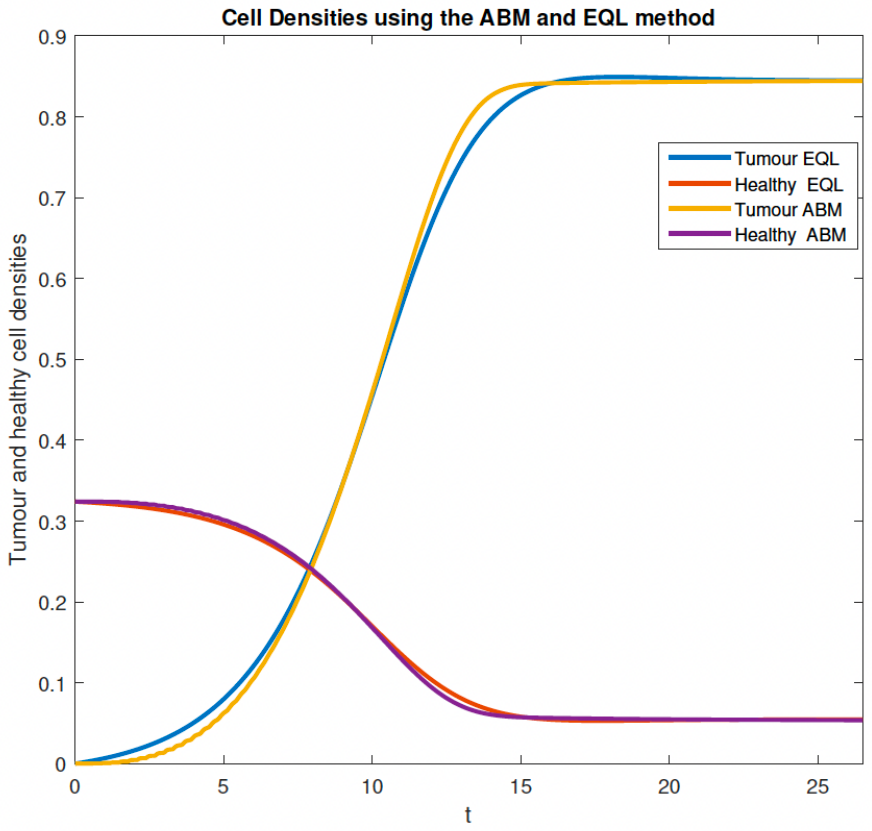}
\label{fig:figEQL}
\caption{Data from tumour and healthy cell densities along with fitted ODE, using MATLAB routine \texttt{lsqr}}
\end{figure}

In the case of the library of chemical reactions, we used all 17 reactions; the non-negative least squares algorithms always worked  and found stable ODE solutions, and determined that the following rate constants were zero
$$ k_1, k_2, k_3, k_4, k_6, k_8, k_{11}, k_{13}, k_{16} $$
when we used either 1800 or 180 time points. The remaining rate constants were

\vspace{0.5cm}
\noindent
1800 time points:
$$ (k_5, k_7, k_9, k_{10}, k_{12}, k_{14}, k_{15}, k_{17}) = (3.2039137, 2.983428, 3.132762, 9.107733, 0.116421,  $$
$$ \quad \quad \quad \quad \quad \quad 0.119287, 3.855036, 4.520207)$$

\noindent
180 time points:
$$ (k_5, k_7, k_9, k_{10}, k_{12}, k_{14}, k_{15}, k_{17}) = (3.063352, 2.586338, 3.000332, 7.912730, 0.113799, $$
$$ \quad \quad \quad \quad \quad \quad 0.079011, 3.44351, 4.45157)$$

The corresponding ODE based on 1800 time points was
\begin{eqnarray}
x' &=& 3.2039x - 3.2491 x^2 + 0.0591 z^2 - 8.3752 x z \nonumber \\
z' &=& 2.9834 z - 9.2296 z^2 + 0.0582 x^2 - 3.8550 x z
\label{eq:eq2star}
\end{eqnarray}
with a Mean-Squared Error (MSE) 
$$ \textrm{MSE} (1.905(-03), 1.996(-04)).$$

The corresponding ODE based on 180 time points was
\begin{eqnarray}
x' &=& 3.0634x - 3.1141 x^2 + 0.0395 z^2 - 7.8951 x z \nonumber \\
z' &=& 2.5863 z - 7.9917 z^2 + 0.0582 x^2 - 3.4435 x z
\label{eq:eq3star}
\end{eqnarray}
with an 
$$ \textrm{MSE} (5.166(-05), 3.278(-06)).$$
The data and the two ODEs based on 1800 time points and 180 time points are given in Figure 5.

\begin{figure}
\centering
$$
\begin{array}{cc}
\includegraphics[trim=0mm 0mm 0mm 0mm,clip,width=0.45\textwidth]{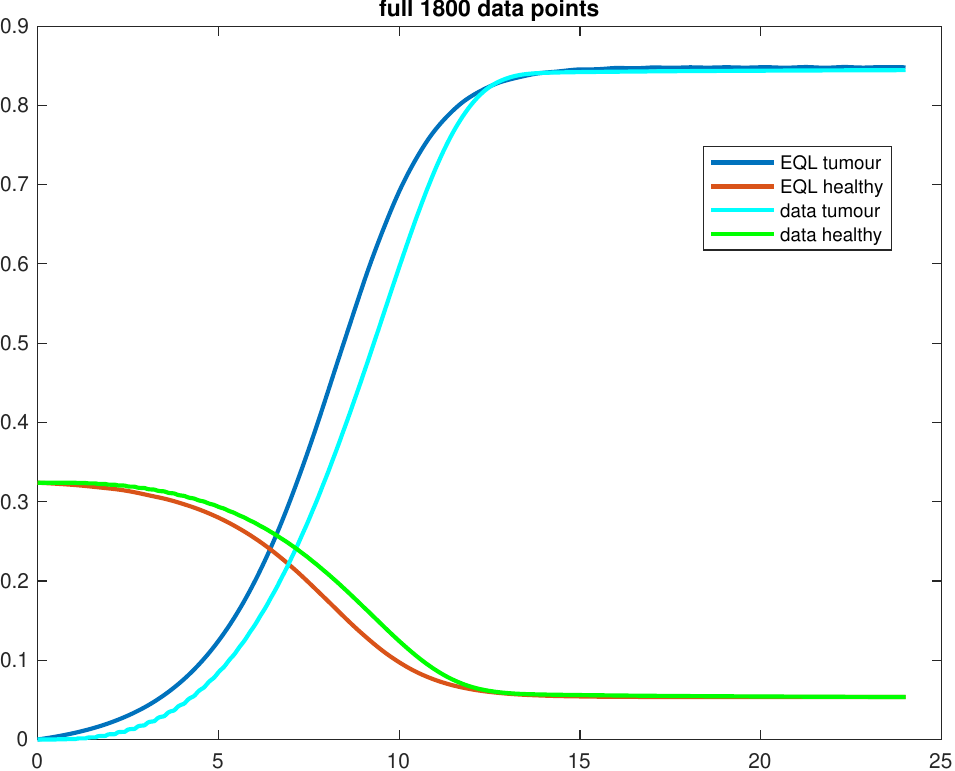} &
\includegraphics[trim=0mm 0mm 0mm 0mm,clip,width=0.45\textwidth]{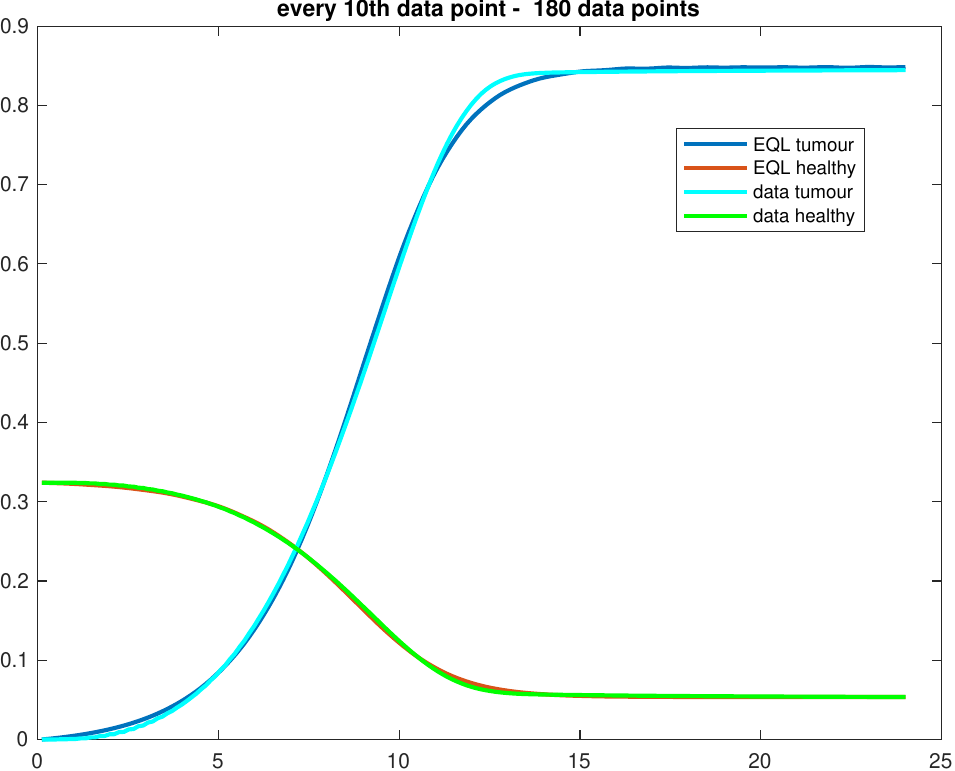} 
\end{array} 
$$
\label{fig:figNewResults}
\caption{Data from tumour and healthy cell densities along with fitted ODE, using the library of chemical reactions and  MATLAB routine \texttt{lsqnonneg}}
\end{figure}

\newpage
\noindent
\textbf{Discussion}

\begin{enumerate}
\item We can see that the solutions obtained from the library of chemical reactions are  very different to the decoupled case. Noticeable differences are that in the coupled case with the library of chemical reactions, there is no $z$ term in the first component and no $x$ term in the second component. Furthermore in the decoupled standard case the $xz$ terms take different signs in the two components.
\item We performed a further exploration in the coupled case to see what would happen if we removed the 12th reaction, namely $X + X \rightarrow Z$ ($k_{12} = 0$). In this case, we still got a very good fit but now the $k_6$ term was non-zero; that is, we had a solution based on the reaction $X \rightarrow X+Z$.  In a further exploration, we removed the 12th and 6th reactions, and now we had the third reaction appearing, associated with a non-zero $k_3$, namely $X \rightarrow Z$. 
\item Finally we removed 12th, 6th and 3rd reactions from the library. In this case there were 11 zero rate constants. In particular $k_2$ and $k_4$ were non-zero, so that the reactions $Z \rightarrow 0, \> Z \rightarrow X$ were now represented in the library. But in this case the MSE was (5.06(-03), 8.036(-4)) which is quite large compared with the other cases, and furthermore the equilibrium values for tumour cells and healthy cells were 0.98 and 1.8(-05), respectively. So this scenario is not a particularly good fit to the data.
\item All three reactions
$$ X+X \rightarrow Z, \quad X \rightarrow X+Z, \quad X \rightarrow Z$$
(while not equivalent) have a similar behaviour and so the data is telling us that a cancer cell should be able to mutate in a healthy cell in some manner.  It is noteworthy that in this case we  do not need a representation from  degradation reactions, as described in (3) above. This is being compensated for by a different type but similar reaction. One reason why this might happen is that the agent-based model does not just consist of a set of chemical reactions, there are other aspects that cannot be represented by simple chemical reactions (for example jumps, and the number of cell divisions for healthy cells). The nonlinear effects seen in Figure 3
also confirm this.
\item Our choice of reactions in the library can potentially lead to issues of linear dependence between reactions -- for example $X + X \rightarrow 0, \> X + X \rightarrow X$. In this case \texttt{lsqnonneg} flags a highly ill-conditioned system, but still converges to a solution in which the reaction $X + X \rightarrow 0$ appears in the final set, but $X + X \rightarrow X$ is removed.
\item We seem to get a better fit to the data when we use less time points -- see Figure 5.
So more data is not necessarily better in terms of the model fit. The great advantage of using a library of chemical reactions is that the non-negative least squares algorithm seems to be able to tell which reactions are important in the data (by finding zero rate constants) and the procedure seems to be very robust in terms of finding a stable ODE system. 
\item One potential criticism of using a library of chemical reactions is that as the number of systems components increases, the size of the library can become large. However, we can use information about how the systems components interact from a reaction viewpoint and this can dramatically reduce the size of the library.
\item To follow up on this final point: we chose a complete library of 17 chemical reactions as we wanted to explore if \texttt{lsgnonneg} would automatically remove reactions that were deemed unnecessary in fitting the data. The solutions found had ony 8 reactions represented in the final ODE system. Considering OCCAM's principle we think this is a nice outcome.
\end{enumerate}

\section{Conclusions}

We have used a simple agent-based model to capture the inherent stochasticity in tumour growth in a fixed two-dimensional spatial domain. This is compute-intensive as many simulations are required to capture, for example, statistics on mean and variance behaviour. However, given a one-off generation of tumour and healthy cell counts we can use equation-learning techniques (SINDy) to derive a differential equation system. We have found  by using a library of chemical reactions in conjunction with a non-negative least squares algorithm that we can robustly obtain a subset of reactions that can characterise the data and in general this leads to a resulting stable ODE. We have also seen that by using different numbers of time points, we can generate a population of close-by differential equations, that closely match the results of the agent-based model. This can allow for uncertainty quantification and an exploration of variability in the model coefficients. 

\vspace{1cm}
\noindent
\textbf{Declarations}

 \noindent
Ethical Approval - not applicable

 \noindent
Availability of supporting data - not applicable

 \noindent
Competing Interests - not applicable

 \noindent
Funding - not applicable

 \noindent
Authors' Contributions - all authors contributed equally.

 \noindent
Acknowledgements - not applicable
 


\end{document}